\documentclass [11pt]{article}
\usepackage{amssymb}
\usepackage{amsmath}
\usepackage{epsfig}
\usepackage{graphicx}
\usepackage{amsthm}
\usepackage{amssymb}
\usepackage{amsfonts}
\usepackage{dirtytalk}
\allowdisplaybreaks
\usepackage{enumerate}
\usepackage{cite}
\theoremstyle{definition}

\begin{document}
\thispagestyle{empty}
\null\vspace{-1cm}
\medskip
\vspace{1.75cm}
\centerline{\textbf{{On Tur\'{a}n Type inequality for Quaternionic Canonical Generalized Polynomials}}}
~~~~~~~~~~~~~~~~~~~~~~~~~~~~~~~~~~~~~~~~~~~~~~~~~~~~~~~~~~~~~~~~~~~~~~~~~~~~~~~~~~~~~~~~~~~~~~~~~~~~~~~~~~~~~~~~~~~~~~~~~~~~~~~~~~~~~~~~~~~~~~~~~~~~~~~~~~~

\centerline{\bf {Idrees Qasim\footnote{Corresponding Author: Idrees Qasim }, Ovaisa Jan}}
\centerline {Department of Mathematics, National Institute of Technology, Srinagar, India-190006}
\centerline {idreesf3@nitsri.ac.in; ovaisa\_2022phamth009@nitsri.ac.in}
~~~~~~~~~~~~~~~~~~~~~~~~~~~~~~~~~~~~~~~~~~~~~~~~~~~~~~~~~~~~~~~~~~~~~~~~~~~~~~~~~~~~~~~~~~~~~~~~~~~~~~~~~~~~~~~~~~~~~~~~~~~~~~~~~~~~~~~~~~~~~~~~~~~~~~~~~
\vskip0.1in
\noindent \textbf{Abstract}: 
This paper establishes Tur\'an-type inequalities for quaternionic canonical generalized polynomials with all zeros in the closed unit ball of radius $k \geq 1$. We extend the classical inequality due to Govil from the complex setting to the quaternionic framework.  We prove that for certain classes of polynomials, the inequality
\[\|P'\| \geq \frac{n}{1+k^n}\|P\|\]
holds for all $k \geq 1$, where $n$ is the degree of the polynomial. Our main results provide sharp derivative estimates for quaternionic canonical generalized polynomials under specific algebraic conditions on their zeros. These findings contribute to the ongoing research of extending classical polynomial inequalities to the noncommutative setting of quaternions.

~~~~~~~~~~~~~~~~~~~~~~~~~~~~~~~~~~~~~~~~~~~~~~~~~~~~~~~~~~~~~~~~~~~~~~~~~~~~~~~~~~~~~~~~~~~~~~~~~~~~~~~~~~~~~~~~~~~~~~~~~~~~~~~~~~~~~~~~~~~~~~~~~~~~~~~~~~

\vspace{0.15in}
~~~~~~~~~~~~~~~~~~~~~~~~~~~~~~~~~~~~~~~~~~~~~~~~~~~~~~~~~~~~~~~~~~~~~~~~~~~~~~~~~~~~~~~~~~~~~~~~~~~~~~~~~~~~~~~~~~~~~~~~~~~~~~~~~~~~~~~~~~~~~~~~~~~~~~~~~

\noindent {{\bf Mathematics Subject Classiﬁcation (2020):} 30G35, 41A17.}\\
\vspace{0.15in}
~~~~~~~~~~~~~~~~~~~~~~~~~~~~~~~~~~~~~~~~~~~~~~~~~~~~~~~~~~~~~~~~~~~~~~~~~~~~~~~~~~~~~~~~~~~~~~~~~~~~~~~~~~~~~~~~~~~~~~~~~~~~~~~~~~~~~~~~~~~~~~~~~~~~~~~~~
\section{Introduction} 

Let \( p(z) = \sum_{v=0}^n a_v z^v \) be a polynomial of degree \( n \), and let \( p'(z) \) denote its derivative.  The extremal problems of functions of complex variables and the results where
some approches to obtaining the classical inequalities are developed on using various
methods of the geometric function theory are known for various norms and for many
classes of functions such as polynomials with various constraints, and on various
regions of the complex plane. Among the most famous results is the Bernstein inequality for the uniform norm on the unit circle. It states:
\[
\max_{|z|=1} |p'(z)| \;\leq\; n  \max_{|z|=1} |p(z)|.
\]
Equality holds if and only if all zeros of \( p \) are at the origin. If the zeros are restricted in some way, the bound can be sharpened. It was conjectured  by Erd\"os and later proved by Lax that if \( p(z) \) has no zeros inside the open unit disk \( |z| < 1 \), then 
\[
\max_{|z|=1} |p'(z)| \;\leq\; \frac{n}{2}  \max_{|z|=1} |p(z)|.
\]
This improves Bernstein's bound by a factor of \( 1/2 \).
In contrast, Turán obtained a lower bound for the maximum of $| p^{\prime}(z)|$ on $|z|=1,$ by proving if \( p(z) \) has all its zeros inside the closed unit disk \( |z| \leq 1 \), then:
\begin{equation}
    \max_{|z|=1} |p'(z)| \;\geq\; \frac{n}{2} \max_{|z|=1} |p(z)|\label{eq:01}
\end{equation}
For both the Erd\"os Lax and Turán inequalities, equality occurs when all zeros lie exactly on the unit circle \( |z| = 1 \).
As an extension of inequality \eqref{eq:01}, Malik \cite{Malik} proved that if $p(z)$ is a polynomial of degree $n$ with complex coefficients having all zeros in $|z|\le k$, $k\le 1$, then
\begin{equation}
	\|p^{\prime}\|\ge \frac{n}{1+k}\|p\|.\label{eq:02}
\end{equation}
\ \\
 Later Govil \cite{Gov} proved, if $p(z)$ is a polynomial of degree $n$ having all zeros in $|z|\le k$, $k\ge 1$, then
\begin{equation}
	\|p^{\prime}\|\ge \frac{n}{1+k^n}\|p\|.\label{eq:03}
\end{equation}
The extension of classical polynomial inequalities from complex analysis to the 
quaternionic setting has emerged as a rich area of research, driven by the 
non-commutative algebra of quaternions and the development of slice regular 
function theory. The foundational work by Gal and Sabadini \cite{Gal} 
established that Bernstein's inequality holds for all quaternionic 
unilateral polynomials of the form $P(q)=\sum q^k a_k$. However, the authors 
demonstrated that the  Erd\"{o}s Lax inequality does not hold in general for quaternionic polynomials. 
They provided the explicit counterexample of polynomial of degree $2$ whose only root is $q=i$ (with multiplicity $2$) and 
which therefore has no zeros in the open unit ball, yet for which $\|P'\| > 
\|P\|$, violating the inequality since $\frac{2}{2}=1$. A partial positive result 
was obtained for polynomials whose zeros are either spheres or real points with 
at most one isolated non-real zero of multiplicity one, a class for which the 
Erd\"{o}s Lax inequality does hold.\\
In \cite{Gal2}, the same authors obtained several partial 
positive results for Tur\'an's inequality under specific algebraic conditions.  Moreover, they observed that when all coefficients of a 
polynomial lie in the same complex plane $\mathbb{C}_I$, Tur\'an's inequality follows directly from the classical 
complex case. Despite these positive results, the authors failed to construct 
counter examples for higher degrees at that time and conjectured that the  Tur\'an's inequality might be valid for all quaternionic polynomials with zeros in the closed unit ball. This conjecture was later disapproved by Coroianu and 
Gal \cite{Coro}. \\
A subsequent comprehensive study by Coroianu and Gal \cite{Coro} 
examined three distinct classes of quaternionic polynomials: canonical 
generalized (QCG) polynomials of the form $(q-\alpha_1)\cdots(q-\alpha_n)$ with 
usual quaternionic multiplication, matrix (QM) polynomials obtained by embedding 
quaternions into $2\times 2$ complex matrices, and unilateral (QU) polynomials 
with the convolution $*$ product. In this paper, for the case of the Bernstein and Erd\"{o}s 
Lax inequalities, results in the larger space of quaternionic matrix polynomials that includes
the space of quaternionic canonical generalized polynomials are obtained. While Tur\'an's inequality  for QCG polynomials was shown to hold for degrees 
$n=1$ and $n=2$ without any additional conditions, and for certain subclasses of 
degree $3$ polynomials when one root is real or two roots coincide, they 
constructed a concrete counterexample for degree $3$ (QCG)-polynomials and proved that the  Turán's inequality cannot be obtained for QCG polynomials in general.
For QU polynomials, the authors proposed several results, subjected to additional assumptions, that establish Turán's inequality. However, for polynomial of degree $3,$ they showed by a counterexample that Turán's inequality for QU polynomials does not hold. Hence, disproved an earlier conjecture, 
raised by Gal and Sabadini \cite{Gal}, that Tur\'an's inequality 
might hold for all quaternionic polynomials with zeros in the closed unit ball.\\
In the paper \cite{LC} Coroianu investigate Bernstein type inequalities for two distinct classes of quaternionic polynomials. First, they show that restrictions of quaternionic unilateral polynomials to any complex plane satisfy a Bernstein inequality for the real Fréchet derivative, which in turn yields a Bernstein inequality for such polynomials under suitable norm conditions. Second, they establish Bernstein type inequalities for quaternionic canonical generalized polynomials, a noncommutative factorized form focusing on cases of degree up to three and on polynomials of arbitrary degree with at most two distinct roots. 
Most recently, Qasim and Jan \cite{OI} proved the  inequality (\ref{eq:02}) and (\ref{eq:03})
to the quaternionic polynomials of degree $n \le 2$ and subclass of quaternionic polynomials of degree $n\ge 3$ for $k\le 1$. The aim of this paper is to extend the  inequality (\ref{eq:02}) and (\ref{eq:03})
to the quaternionic polynomials of degree $n \le 2$ and subclass of quaternionic polynomials of degree $n \ge 3$ for $k\ge 1.$
\section{Preliminary}
We begin with some preliminaries on quaternions and quaternionic polynomials. Quaternions, a 4-dimensional division algebra discovered by Hamilton in 1843, are discussed herein along with their definitions and significant properties. Quaternions are typically denoted by $\mathbb{H}.$\\
\textbf{Quaternions}:\label{def1.2.1}
     A quaternion \(q \in \mathbb{H}\) is an expression of the form

\[
q = q_0 + q_1 {i} + q_2 {j} + q_3 {k}
\]
where \(q_0, q_1, q_2, q_3 \in \mathbb{R}\) (real numbers). The imaginary units satisfy

\[
{i}^2 = {j}^2 = {k}^2 = {i}{j}{k} = -1
\]
A quaternion can also be represented as:
$q= q_{0}+ v,  ~\mbox{where}~ q_{0}\in \mathbb{R} ~\mbox{and}~ v = q_{1} \hat{i} + q_{2}\hat{j}+q_{3}\hat{k}$.\\
\textbf{Conjugation}:{\label{def1.2.2}}
The quaternion conjugate of \(q \in \mathbb{H}\) is given by

\[
q^{*} = q_0 - q_1 {i} - q_2 {j} - q_3 {k}. 
\]\\
    \textbf{Norm of a quaternion}:{\label{def1.2.3}} The norm of a quaternion \(q \in \mathbb{H}\) is given by the real number

\[
|q|^2 = q {q^{*}} = q_0^2 + q_1^2 + q_2^2 + q_3^2
\]
The norm satisfies \(|pq| = |p||q|\) and is multiplicative. \\
\textbf{Quaternionic polynomial}:{\label{def1.2.4}}
     A quaternionic polynomial of degree $n$ is the expression $f(q)=\sum_{k=0}^{n}a_kq^k$ or $P(q)=\sum_{k=0}^{n}q^ka_k,$ $a_n\neq 0$, $a_k\in \mathbb{H},$ $k=0,1,\dots,n$ in the quaternionic indeterminate $q$. These are commonly termed as simple quaternionic polynomials. Unlike complex case, the multiplication in $\mathbb{H}$ is non-commutative, one can consider polynomials with coefficients on the left or on the right. Polynomials with coefficients on one side are called quaternionic unilateral polynomials. The derivative of a polynomial $P(q)=\sum\limits_{k=0}^{n}q^{k}a_k$ is $P^{\prime}(q)=\sum\limits_{k=1}^{n}kq^{k-1}a_k.$ For more details on quaternionic polynomials see \cite{Lam, LC, Ab3}. Two quaternionic polynomials of the type $P(q)=\sum_{k=0}^{n}q^ka_k$ can be multiplied according to the convolution product (Cauchy multiplication rule) as given $P_1(q)=\sum\limits_{i=0}^{n}q^ia_i$ and $P_2(q)=\sum\limits_{j=0}^{m}q^jb_j$, we define
$$(P_1*P_2)(q):=\sum_{i=0,1,\dots, n~~j=0,1,\dots,m}^{}q^{i+j}a_ib_j.$$
If $P_1$ has real coefficients, then so called * multiplication coincides with the usual pointwise multiplication.\\
\indent The lack of commutativity results in a behavior of polynomials that is quite distinct from their behavior in the real or complex settings. For instance, a real or complex polynomial of degree $n$ can have at most $n$ (real or complex) zeros, counted with their multiplicity. However, in the quaternionic setting, the second-degree polynomial $q^2+1$ exhibits an infinite number of zeros. We summarize several results on the zeros of quaternionic polynomials in the following theorem, see \cite{Gor, Lam}.\\\\
\textbf{Theorem 2.1} \begin{enumerate}
	\item Let $P(q)$ be a non-zero polynomial. A quaternion $\alpha$ is a zero of $P(q)$ if and only if $q-\alpha$ is a left divisor of $P(q)$.
	\item If $P(q)=(q-\alpha_1)*\dots *(q-\alpha_n)a_n$, where $\alpha_i\in \mathbb{H},~i=1,2,\dots, n$, then $\alpha_1$ is the zero of $P$ and every other zero of $P$ belongs to the spheres $[\alpha_i],~i=1,2,\dots,n$.
	\item If $P(q)$ is a non-zero polynomial of degree $n$, then there exists $\alpha_1,\dots,\alpha_n$ such that $P(q)=(q-\alpha_1)*\dots*(q-\alpha_n)a_n$
\end{enumerate}
\textbf{2-Dimensional sphere}:{\label{def1.2.31}} Let $$\mathbb{S}:=\left\{q=q_{0}+q_{1}i+q_{2}j+q_{3}k\in\mathbb{H}:q_{0}=0, ~q_{1}^{2}+q_{2}^{2}+q_{3}^{2}=1\right\}$$ be the unit sphere of purely imaginary quaternions. The set $\mathbb{S}$ is a 2-dimensional sphere in $\mathbb{H}$ identified with $\mathbb{R}^{4}$. The elements $I\in\mathbb{S}$ are called imaginary units as $I^{2}=-1.$ In particular, $i=(0,1,0,0),~j=(0,0,1,0),~k=(0,0,0,1)\in \mathbb{S}$.
\noindent For some $x,y\in \mathbb{R},~y>0$ and some $I\in\mathbb{S},$ every quaternion $q\notin \mathbb{R}$ can be uniquely expressed as $q=x+yI.$
In other words (denoting $\mathbb{R} = \mathbb{R} \cdot 1 \subset \mathbb{H}$), the quaternion algebra $\mathbb{H}$ can be written as the union
\[
\mathbb{H} = \bigcup_{I \in \mathbb{S}} (\mathbb{R} + \mathbb{R}I)
\]
of complex planes $L_I := \mathbb{R} + \mathbb{R}I$, each isomorphic to $\mathbb{C}$, and all sharing only the real axis. This decomposition neatly reflects the non-commutativity of $\mathbb{H}$: two quaternions commute if and only if they lie in the same complex plane $L_I;$ which is always true when one of them is real. \\
Let $\Omega$ be a domain in the space of quaternions $\mathbb{H}$ and let
$\mathbb{S}= \left\{q\in \mathbb{H} : q^{2}= - 1\right\}$
denote the 2-sphere of quaternionic imaginary units. We define the notion of regular
functions as follows:\\
\textbf{ Slice Regular function}:{\label{def1.2.6}}
    Let $f$ be a quaternion-valued function defined on a slice domain $\Omega \subseteq \mathbb{H}$. For an imaginary unit $I\in \mathbb{S}$, let $\Omega_{I}=\Omega\cup C_{I}$, where $C_{I}:=\mathbb{R}+\mathbb{R}I$, and let $f_{I}=f|_{\Omega_{I}}$ be the restriction of $f$ to $\Omega_{I}$. The restriction $f_{I}$ is called holomorphic if it has continuous partial derivatives and 
$$\frac{1}{2}\left(\frac{\partial}{\partial x}+I\frac{\partial}{\partial y}\right)f_{I}(x+Iy)\equiv 0.$$
The function $f$ is called slice regular if for all $I\in \mathbb{S}$, $f_{I}$ is holomorphic.\\
\noindent The regular functions of a quaternionic variable have been intensively studied in the past decade and their rapid development has been largely driven by the application to operator theory.\\
From now onwards, for a quaternionic polynomial $P$, the norm of $P$ is defined by
$$\|P\|=\max \{|P(q)|:|q|\le 1\}=\max \{|P(q)|:|q|=1\}.$$

	\section{Main Results}
In this section, we prove inequality \eqref{eq:03} for subclass of quaternionic polynomials for $n\ge 3$. In fact, we prove:\\
\noindent \textbf{Theorem 3.1.}
	Let 
	\[
	P(q)=(q-\alpha)(q-\beta)(q-\gamma)
	\]
	be a quaternionic polynomial of degree $3$, where
	\[
	|\alpha|\le k,\qquad |\beta|\le k,\qquad |\gamma|\le k,
	\]
	If 
	\[
	k \ge 2^{1/3},
	\]
	then
	\[
	\|P'\| \;\ge\; \frac{3}{1+k^{3}}\,\|P\|.
	\]

\begin{proof}
	Set
	\[
	B=\alpha\beta+\alpha\gamma+\beta\gamma,
	\qquad 
	C=\alpha\beta\gamma.
	\]
	Expansion gives
	\[
	P(q)=q^{3}+qB-C,
	\qquad 
	P'(q)=3q^{2}+B.
	\]
	For $|q|=1$,
	\[
	|P(q)|\le |q^{3}|+|qB|+|C|
	\le 1+|B|+|C|,
	\]
	hence
	\begin{equation}
		\|P\|\le 1+|B|+|C|.
		\label{eq:Pupper}
	\end{equation}
	If $B=0$, then $P'(q)=3q^{2}$, so $\|P'\|=3$.
	If $B\neq 0$, pick $q_{0}$ on the unit sphere with $q_{0}^{2}=B/|B|$.  
	Then
	\[
	P'(q_{0}) = 3\frac{B}{|B|} + B 
	= (3+|B|)\frac{B}{|B|},
	\]
	so
	\[
	|P'(q_{0})| = 3+|B|.
	\]
	As $|P'(q)|\le 3+|B|$ for all $|q|=1$, we obtain
	\begin{equation}
		\|P'\|=
		\begin{cases}
			3, & B=0,\\
			3+|B|, & B\ne 0.
		\end{cases}
		\label{eq:Pprime}
	\end{equation}
	\medskip
	\textbf{ Case $B=0$.}
	From \eqref{eq:Pupper} and $|C|\le k^{3}$,
	\[
	\|P\|\le 1+k^{3}.
	\]
	Therefore,
	\[
	\frac{3}{1+k^{3}}\|P\|
	\;\le\;
	\frac{3}{1+k^{3}}(1+k^{3})
	=3 = \|P'\|.
	\]
	\medskip
	\textbf{Case $B\ne 0$.}
	Using \eqref{eq:Pupper} and \eqref{eq:Pprime}, it is sufficient to prove
	\begin{equation}
		3+|B|
		\;\ge\;
		\frac{3}{1+k^{3}}(1+|B|+|C|).
		\label{eq:keyineq}
	\end{equation}
	Since $|C|\le k^{3}$, a stronger inequality is
	\[
	3+|B|\;\ge\;\frac{3}{1+k^{3}}(1+|B|+k^{3}),
	\]
	i.e.
	\[
	(3+|B|)(1+k^{3})
	\;\ge\;
	3(1+|B|+k^{3}).
	\]
	Expanding both sides gives
	\[
	3k^{3}+|B|(k^{3}-2)\ge 0.
	\]
	If $k\ge 2^{1/3}$, then $k^{3}\ge 2$, so $k^{3}-2\ge 0$ and the above inequality holds for all $|B|\ge 0$.  
	Therefore \eqref{eq:keyineq} holds whenever $k\ge 2^{1/3}$.\\
	\medskip
	Thus, for every $k\ge 2^{1/3}$,
	\[
	\|P'\| = 3+|B|\;\ge\; \frac{3}{1+k^{3}}\|P\|,
	\]
	completing the proof.
\end{proof}

\noindent \textbf{Theorem 3.2.} Let
	$$P(q)=(q-\alpha)*(q-\beta)*(q-\gamma)$$
	be a quaternionic polynomial of degree $3$, where
	$$
	|\alpha|\le k,\qquad |\beta|\le k,\qquad |\gamma|\le k.
	$$ If $\alpha+\beta+\gamma= 0$ and \begin{equation}
	\left\{|\alpha|\le A(k)~or~|\beta|\le A(k)~or~|\gamma|\le A(k)\right\},\label{eq:901}
\end{equation}
    Then for $k \ge 1,$
	
	$$
	\|P'\| \;\ge\; \frac{3}{1+k^{3}}\,\|P\|.
	$$
    where $A(k)=
\dfrac{-3k^{2}+\sqrt{\,9k^{4}-4(2-k^{3})\bigl((2-k^{3})k^{2}-3k^{3}\bigr)\,}}
{2(2-k^{3})}.$ 
    \begin{proof} Expanding $P$ and using $\alpha +\beta +\gamma =0 $ gives\\ $$P(q)= q^{3}+qB-C;\quad P^{\prime}(q)=3q^{2}+B$$ 
	where,
	\[
	B=\alpha\beta+\alpha\gamma+\beta\gamma,
	\qquad 
	C=\alpha\beta\gamma.
	\]
For $|q|=1$,
	\[
	|P(q)|\le |q^{3}|+|qB|+|C|
	\le 1+|B|+|C|,
	\]
	hence
	\begin{equation}
		\|P\|\le 1+|B|+|C|.
		\label{eq:04}
	\end{equation}
	If $B=0$, then $P'(q)=3q^{2}$, so $$\|P'\|=3$$
 If $B\neq 0$, pick $q_{0}$ on the unit sphere with $q_{0}^{2}=B/|B|$.  
	Then
	\[
	P'(q_{0}) = 3\frac{B}{|B|} + B 
	= (3+|B|)\frac{B}{|B|},
	\]
	so
	\[
	|P'(q_{0})| = 3+|B|.
	\]
	As $|P'(q)|\le 3+|B|$ for all $|q|=1$, we obtain
	\begin{equation}
		\|P'\|=
		\begin{cases}
			3, & B=0,\\
			3+|B|, & B\ne 0.
		\end{cases}
		\label{eq:05}
	\end{equation}
	\textbf{ Case $B=0$.}
	From \eqref{eq:04} and $|C|\le k^{3}$,
	\[
	\|P\|\le 1+k^{3}.
	\]
	Therefore,
	\[
	\frac{3}{1+k^{3}}\|P\|
	\;\le\;
	\frac{3}{1+k^{3}}(1+k^{3})
	=3 = \|P'\|.
	\] So inequality holds.\\
    
	\noindent \textbf{ Case $B\ne 0$.}
	Using \eqref{eq:04} and \eqref{eq:05}, it is sufficient to prove
	\begin{equation*}
		3+|B|
		\;\ge\;
		\frac{3}{1+k^{3}}(1+|B|+|C|).
	\end{equation*}
	This is equivalent to $$(1+k^{3})(3+|B|)\ge 3(1+|B|+|C|)$$
    which simplifies to \begin{equation}
        3k^{3}+(k^{3}-2)|B|\ge 3|C|
    \label{eq:06} \end{equation}  Using $\alpha+\beta+\gamma=0$ gives $\beta+\gamma=-\alpha,$ hence $$B=\alpha(\beta+\gamma)+\beta\gamma=-\alpha^{2}+\beta \gamma$$ Applying the triangle inequality $$|B|\le |\alpha|^{2}+|\beta \gamma|\le|\alpha|^{2}+k^{2}
    $$ and $$|C|=|\alpha\beta\gamma|\le|\alpha|k^{2}$$ Now inequality (\ref{eq:06}) will become $$3k^{3}+(k^{3}-2)(|\alpha|^{2}+k^{2})\ge 3|\alpha|k^{2}$$ A sufficient condition for (\ref{eq:06}) is therefore \begin{equation}
        (k^{3}-2)x^{2}-3k^{2}x+k^{2}(k^{3}-2)+3k^{3}\ge0,~ x=|\alpha| 
    \end{equation}{\label{eq:07}}
This simplifies to \begin{equation}
    (2-k^{3})x^{2}+3k^{2}x+(2-k^{3})k^{2}-3k^{3}\le 0 \label{eq:08}
\end{equation}
Define the quadratic function $$\phi_{k}(x)=(2-k^{3})x^{2}+3k^{2}x+(2-k^{3})k^{2}-3k^{3}$$
\textbf{Subcase I:} If $1\le k < 2^{\frac{1}{3}}$ then $2-k^{3}>0$. Hence $\phi_{k}$ is an upward opening parabola.\\
Since $\phi_{k}(0)=(2-k^{3})k^{2}-3k^{3}<0,$ for $k\ge 1$ there exists a positive root. The positive root is given by \\
$$A(k)= \dfrac{-3k^{2}+\sqrt{9k^{4}-4(2-k^{3})[(2-k^{3})k^{2}-3k^{3}]}}{2(2-k^{3})}$$
Then for $|\alpha|\le A(k),$ we have $\phi_{k}(x)\le 0.$ So (\ref{eq:08}) holds implying (\ref{eq:06}) and thus the desired inequality. By symmetry the same conclusion holds if $|\beta|\le A(k)$ or $|\gamma|\le A(k).$\\

\noindent \textbf{Subcase II:} If $K = 2^{\frac{1}{3}}$ then $2-k^{3}=0$ Thus (\ref{eq:08}) reduces to $3k^{2}x\le 3k^{3}$ i.e $x\le k$ which is true. Since $|\alpha|\le k.$ Hence (\ref{eq:08}) holds automatically and the desired inequality follows.\\
\noindent\textbf{Subcase III:} If $k\ge 2^{\frac{1}{3}}$ then $ 2-k^{3}<0.$ Here $2-k^{3}<0.$ Thus $\phi_{k}$ is a downward opening parabola $$\phi_{k}(0)= (2-k^{3})k^{2}-k^{3}<0$$ holds for all $x\in [0,k].$ Hence the inequality holds.
	Thus, for every $k\ge 1$,
	\[
	\|P'\| = 3+|B|\;\ge\; \frac{3}{1+k^{3}}\|P\|,
	\]
\end{proof}

\noindent \textbf{Theorem 3.3.} Let 
	\[
	P(q)=(q-\alpha)*(q-\beta)*(q-\gamma)*(q-\lambda)
	\]
	be a quaternionic polynomial of degree \(4\), where
	\[
	|\alpha|,\,|\beta|,\,|\gamma|,\,|\lambda| \le k,\qquad k\ge 1.
	\]
	Assume that the first two symmetric sums vanish:
	\[
	\alpha+\beta+\gamma+\lambda = 0, \qquad
	\alpha\beta+\alpha\gamma+\alpha\lambda+\beta\gamma+\beta\lambda+\gamma\lambda = 0.
	\]
	Set 
	\[
	E=\alpha\beta\gamma+\alpha\beta\lambda+\alpha\gamma\lambda+\beta\gamma\lambda,
	\qquad 
	C=\alpha\beta\gamma\lambda.
	\]
	Then:
	
	\begin{enumerate}[(i)]
		\item If \(1 \le k \le 3^{1/4}\), define
		\[
		x^* = \bigl(A+B\bigr)^{1/3} + \bigl(A-B\bigr)^{1/3},
		\]
		where
		\[
		A = -\,\frac{(3-k^4)k^3 - 4k^4}{2(3-k^4)}, 
		\quad
		B = \sqrt{\frac{\bigl((3-k^4)k^3 - 4k^4\bigr)^2}{4(3-k^4)^2} 
			+ \frac{(4k^3)^3}{27(3-k^4)^3}}.
		\]
		If one of the following holds:
		\[
		|\alpha|\le x^*,\qquad 
		|\beta|\le x^*\ (\beta\in\mathbb R),\qquad
		|\gamma|\le x^*\ (\gamma\in\mathbb R),\qquad
		|\lambda|\le x^*,
		\]
		then
		\[
		\|P'(q)\|\;\ge\;\frac{4}{\,1+k^4\,}\,\|P(q)\|.
		\]
		
		\item If \(k \ge 3^{1/4}\), then the inequality
		\[
		\|P'(q)\|\;\ge\;\frac{4}{\,1+k^4\,}\,\|P(q)\|
		\]
		holds unconditionally (i.e.\ no small–root assumption is needed).
	\end{enumerate}
\begin{proof}
	Using the hypotheses on the first two symmetric sums, 
	\[
	P(q)=q^4 - qE + C,
	\qquad
	P'(q)=4q^3 - E,
	\]
	where \(E\) and \(C\) are as above.  
	For \(|q|=1\), 
	\[
	|P(q)| \le 1 + |E| + |C|,
	\qquad
	|P'(q)| \le 4 + |E|.
	\]
	
	\smallskip\noindent
	Choose \(q_0\) on the unit sphere such that 
	\[
	q_0^3 = -\,\frac{E}{|E|}.
	\]
	Then
	\[
	P'(q_0)=4q_0^3 - E 
	= -4\frac{E}{|E|} - E
	= -(4+|E|)\frac{E}{|E|},
	\]
	and therefore
	\begin{equation}
	    \|P'\| = 4 + |E|.
	\label{eq:09}
	\end{equation}

	\smallskip\noindent
	Since \(\|P\|\le 1 + |E| + |C|\), it suffices to ensure
	\[
	4+|E|
	\;\ge\;
	\frac{4}{1+k^4}\bigl(1+|E|+|C|\bigr).
	\]
	Multiplying by \(1+k^4\) gives the equivalent condition
	\begin{equation}
	    4k^4 
	\;\ge\; 
	(3-k^4)|E| + 4|C|.
	\label{eq:10}
	\end{equation}
	
	\smallskip\noindent
	\textbf{Bounding $E$ and $C$ using the symmetric identities.}
	From \[
	\alpha+\beta+\gamma+\lambda = 0, \qquad
	\alpha\beta+\alpha\gamma+\alpha\lambda+\beta\gamma+\beta\lambda+\gamma\lambda = 0,
	\] one checks (as in the $k\le1$ case) that
	\[
	E = \alpha^3 + \beta\gamma\lambda,
	\qquad 
	|C| = |\alpha\beta\gamma\lambda| \le |\alpha|\,k^3,
	\qquad
	|\beta\gamma\lambda| \le k^3.
	\]
	Thus
	\begin{equation}
	    (3-k^4)|E| + 4|C|
	\;\le\; (3-k^4)|\alpha|^3 + (3-k^4)k^3 + 4|\alpha|k^3.
	\label{eq:11}
	\end{equation}
	
	\smallskip\noindent
	\textbf{Case 1: \(1 \le k \le 3^{1/4}\).}
	Here \(3-k^4 \ge 0\).  
	Define
	\[
	S_k(x) := (3-k^4)x^3 + 4k^3x + (3-k^4)k^3 - 4k^4.
	\]
	Then by \eqref{eq:11}, inequality \eqref{eq:10} is guaranteed whenever
	\[
	S_k(|\alpha|) \le 0.
	\]
	Because
	\[
	S'_k(x) = 3(3-k^4)x^2 + 4k^3 > 0,
	\]
	the function \(S_k\) is strictly increasing on \([0,\infty)\); hence there is a unique positive root \(x^*\), and
	\[
	S_k(x)\le0 \iff 0\le x\le x^*.
	\]
	Applying Cardano’s formula gives the stated expression for \(x^*\).  
	Therefore, if \(|\alpha|\le x^*\), inequality \eqref{eq:10} holds.  
	By symmetry, the same conclusion holds if \(|\beta|\le x^*\), \(|\gamma|\le x^*\), or \(|\lambda|\le x^*\), yielding
	\[
	\|P'\| \;\ge\; \frac{4}{1+k^4}\,\|P\|.
	\]
	
	\smallskip\noindent
	\textbf{Case 2: \(k \ge 3^{1/4}\).}
	Here \(k^4 - 3 \ge 0\).  
	Using \eqref{eq:09},
	\[
	4+|E| = \|P'\|,\qquad |C|\le k^4.
	\]
	Thus  
	\[
	4+|E|
	\;\ge\;
	\frac{4}{1+k^4}(1+|E|+k^4)
	\quad\Longleftrightarrow\quad
	(k^4-3)|E| \ge 0,
	\]
	which is automatically true when \(k\ge 3^{1/4}\).  
	Hence the desired estimate holds with no additional assumptions in this range.
	\smallskip
	This completes the proof.
\end{proof}

	\section{\large{Declaration}}
\textbf{Availabilty of data and material}\\
 Data availability is not applicable to this article as no new data were created or analyzed in this study.\\

\noindent \textbf{Competing interests}\\
The author declare that they have no competing interests.\\

\noindent \textbf{Funding}\\
The research of second author is supported by DST Inspire Fellowship (ID No. IF210629).

\end{document}